\documentclass[a4paper,12pt]{amsart}
\usepackage{pst-node}
\usepackage{subfigure}

\newgray{lightgray}{0.9}  

\theoremstyle{plain}

\newtheorem{lemma}{Lemma}
\newtheorem*{corollary}{Corollary}
\newtheorem*{main theorem}{Theorem}

\theoremstyle{definition}
\newtheorem*{acknowledgements}{Acknowledgements}

\newtheorem{remark}{Remark}
\newtheorem{example}{Example}
\newtheorem*{notation}{Notation}

\DeclareMathOperator{\T}{T}

\newcommand{\field}[1]{\mathbb{#1}}
\newcommand{\C}{\field{C}}
\newcommand{\R}{\field{R}}
\newcommand{\Z}{\field{Z}}
\newcommand{\N}{\field{N}}

\begin{document}

\title{Frankel's theorem in the symplectic category}
\author{Min Kyu Kim}
\address{Department of Mathematics,
Korea Advanced Institute of Science and Technology, 373-1,
Kusong-Dong, Yusong-Gu, Taejon, 305-701, Korea}
\email{minkyu@kaist.ac.kr}



\begin{abstract}
We prove that if an $(n-1)$-dimensional torus acts symplectically
on a $2n$-dimensional symplectic manifold, then the action has a
fixed point if and only if the action is Hamiltonian. One may
regard it as a symplectic version of Frankel's theorem which says
that a K\"{a}hler circle action has a fixed point if and only if
it is Hamiltonian. The case of $n=2$ is the well known theorem by
McDuff.
\end{abstract}

\maketitle

\section{Introduction} \label{section: 1}

One of interesting problems in symplectic geometry is to find
conditions to guarantee that a symplectic action must be
Hamiltonian. It has a long history since the beautiful theorem of
Frankel which says that a K\"{a}hler circle action has a fixed
point if and only if it is Hamiltonian \cite{Fr}, \cite{O}. For
the history of the problem, see \cite{M}, \cite[Section 1, Chapter
5]{MS}. In 1988, McDuff constructs a six-dimensional symplectic
non-Hamiltonian circle action with fixed tori \cite{M}, and it
shows that Frankel's theorem is not true in the symplectic
category. So, one may think that we need more symmetry for
existence of a fixed point to guarantee that a symplectic action
must be Hamiltonian. In this report, we prove the following.

\begin{main theorem}
Let an $(n-1)$-dimensional torus act symplectically on a
$2n$-dimensional symplectic manifold in an effective way. Then the
action has a fixed point if and only if the action is Hamiltonian.
\end{main theorem}

From the McDuff's example and its products with copies of a
two-dimensional sphere endowed with the usual rotation, we can see
that the condition on the dimension of the acting torus is optimal
to obtain the theorem. For $(n-1)$-dimensional torus actions on
$2n$-dimensional symplectic manifolds, see \cite{KT1}, \cite{KT2}.

To guarantee that a symplectic action must be Hamiltonian, we need
some condition either on the manifold or on the action as noted in
\cite[page 155]{MS}. For the former, Feldman shows that a
symplectic circle action with nonempty fixed points on a manifold
with the positive Todd genus is Hamiltonian \cite{Fe}. For the
latter, Tolman and Weitsman show that a semifree symplectic circle
action with nonempty isolated fixed points is Hamiltonian
\cite{TW}. It is also conceivable that a symplectic circle action
with nonempty isolated fixed points must be Hamiltonian, but it is
still open. In this report, we impose a restriction on the acting
group itself. Recently, Sleewaegen reproved the theorem of this
report under an additional assumption in \cite{Sl}.

In Section~\ref{section: 2}, we explain notations and define
terminologies. Also, local behavior of generalized moment maps is
investigated. The proof is given in Section~\ref{section: 3}. It
is based on analysis of local behavior of generalized moment maps.
Also, McDuff's note \cite[proof of Lemma 2]{M} is repeatedly used
which says that a symplectic circle action with a local extremum
for a generalized moment map must be Hamiltonian.

\begin{acknowledgements}
I would like to give special thanks to Dong Youp Suh for bringing
the paper of Tolman and Weitsman to my attention. It was the
beginning of this research.

I would like to thank Jin-Whan Yim for many interesting
discussions. I also would like to thank Seung-Taik Oh and
Myung-Jun Choi for teaching me how to draw pictures.

\end{acknowledgements}

\section{Local description of generalized moment maps} \label{section: 2}
It is well known that moment map images classify toric manifolds
\cite{D}. Similarly, in her beautiful paper \cite{T} Tolman shows
that moment map images(she calls them x-rays) are very useful in
dealing with six-dimensional Hamiltonian two-torus actions. In
this section, we define the x-ray of an action and investigate
local behavior of generalized moment maps through observing
x-rays.


 First, we explain notations. Let $(M^{2n}, \omega)$ be a
$2n$-dimensional symplectic manifold. Let us fix a decomposition
of the $(n-1)$-torus $T^{n-1}=S^1_1 \times \cdots \times
S^1_{n-1}$ where $S^1_i,$ $i=1, \cdots, n-1,$ are circle subgroups
of $T^{n-1}.$ We denote by $\mathfrak t$ the Lie algebra of
$T^{n-1}.$ Let $\Lambda_0$ be the kernel of the exponential map
for $T^{n-1},$ i.e. the lattice of circle subgroups of $T^{n-1}.$
We assume that $T^{n-1}$ acts symplectically on $M^{2n}$ in an
effective way with a nonempty fixed point set. We may assume that
the symplectic form $\omega$ is integral, and hence admits a
generalized moment map \cite[Lemma 1]{M}. Let $\R/\Z$ valued
functions $\mu_i,$ $i=1, \cdots, n-1,$ be generalized moment maps
for the $S^1_i$ actions on $M^{2n}.$ We denote the range of
generalized moment maps $\mu_i$ by $\R/\Z$ instead of $S^1$ to
avoid confusion with groups acting on manifolds. Put $\mu=(\mu_1,
\cdots, \mu_{n-1}).$ Locally, the range $(\R/\Z)^{n-1}$ of the
generalized moment map $\mu$ can be considered as ${\mathfrak
t}^*,$ and so every local result on moment maps including Local
Convexity Theorem \cite{GS} is also true for generalized moment
maps. For each $x$ in the fixed point set $M^{T^{n-1}}$, we denote
by $\T_x M$ the tangent space of $M^{2n}$ at $x.$ We also denote
by the same notation $\T_x M$ the linear isotropy representation
of $T^{n-1}$ on the tangent space at $x,$ and by $\alpha_{i,x},$
$i=1, \cdots, n,$ the weights of $\T_x M.$ The weights
$\alpha_{i,x}$ can be regarded as elements of $\mathfrak t^*$
through differentiation.

For a symplectic circle action, it is known that the action is
non-Hamiltonian if and only if each component of the fixed point
set is not a local extremum for a generalized moment map of the
action \cite[proof of Lemma 2]{M}, and this is equivalent to say
that for each fixed point $x$ the cone generated by the weights
$\alpha_{i,x}$ is equal to the whole $(L S^1)^*$ where $L S^1$
means the Lie algebra of $S^1.$ Similarly, for a $2n$-dimensional
symplectic $(n-1)$-torus action, the cone generated by weights of
the linear isotropy representation at each fixed point plays a key
role in the report, and so the following notations are used.
\begin{notation}
For $\alpha_i \in {\mathfrak t}^* , i=1, \cdots , r ,$ let
$$S(\alpha_1 , \cdots, \alpha_r )=\{s_1
\alpha_1 +\cdots+s_r \alpha_r \in \mathfrak t^* | s_1, \cdots, s_r
\ge 0 \} ,$$
$$S^\circ (\alpha_1 , \cdots , \alpha_r )=\{ s_1 \alpha_1 + \cdots
+ s_r \alpha_r \in {\mathfrak t}^* | s_1 , \cdots , s_r >0  \} .$$
\end{notation}

Now, we define the x-ray of an action. Let $T_1 , \cdots , T_N$ be
the subgroups of $T^{n-1}$ which occur as stabilizers of points in
$M^{2n}$. Let $M_i$ be the set of points for which the stabilizer
is $T_i.$ We also denote such a set by $M_{T_i}.$ By relabeling we
can assume that $M_i$'s are connected and the stabilizer of points
in $M_i$ is $T_i.$ Then, $M^{2n}$ is a disjoint union of $M_i$'s.
Also, it is well known that $M_i$ is open dense in its closure and
the closure is just a component of $M^{T_i}.$ Let $\mathfrak{M}$
be the set of $M_i$'s. Then, the \textit{x-ray} of $(M^{2n},
\omega, \mu)$ is defined as the set of $\mu(\overline{M_i})$'s for
$M_i \in \mathfrak{M}.$ Each image $\mu(\overline{M_i})$(resp.
$\mu(M_i)$) is called an \textit{$m$-face}(resp. \textit{an open
$m$-face}) of the x-ray if $T_i$ is $(n-1-m)$-dimensional. Our
interest is in open $(n-2)$-faces of the x-ray which are of
codimension one in $(\R/\Z)^{n-1}$ by \cite[Theorem 3.6]{GS}. See
also Example~\ref{example: simple form of moment} below.

The generalized moment map has a simple form in a neighborhood of
each orbit. Let $x$ be a point of $M$ with the stabilizer
$T_x^{n-1}.$ The \textit{symplectic slice representation} $V$ at
$x$ is defined as the induced $T_x^{n-1}$ representation
$$(\T_x T^{n-1} \cdot x)^\omega / \T_x (T^{n-1} \cdot x)$$ where
the superscript $\omega$ means the symplectic perpendicular with
respect to the form $\omega.$ Then, a neighborhood of the orbit
$T^{n-1} \cdot x$ is equivariantly symplectomorphic to a
neighborhood of
$$E = T^{n-1} \times_{T_x^{n-1}} ( \mathfrak{t}_x^\circ \times V)$$
where $\mathfrak{t}_x^\circ$ is the annihilator of
$\mathfrak{t}_x$ in $\mathfrak{t}^*$ and $T_x^{n-1}$ acts
trivially on $\mathfrak{t}_x^\circ.$ In the neighborhood, the
generalized moment map is given by
$$\mu([t, \eta, v])=\mu(x)+ \eta+ A^* \mu_V (v)
$$ where $A: \mathfrak{t} \rightarrow \mathfrak{t}_x$ is a
projection, $A^*$ is dual to $A,$  and $\mu_V : V \rightarrow
\mathfrak{t}_x^*$ is the moment map for the symplectic slice
representation $V.$ For more details, see \cite[p.
4209--4210]{LT}.

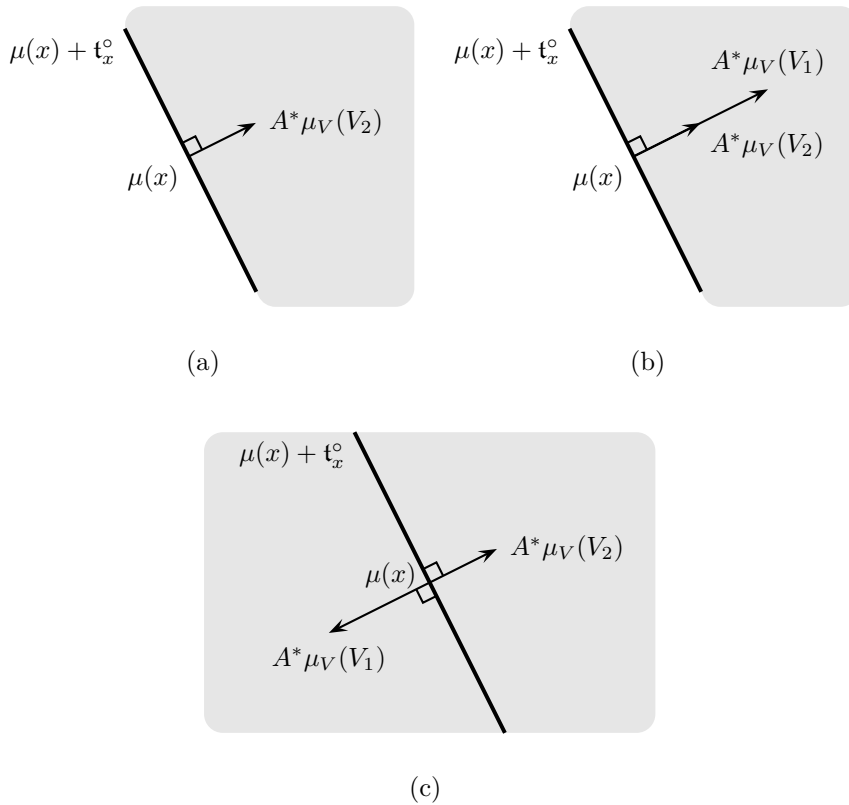
\begin{figure}[ht]
\begin{center}
\mbox{

\subfigure[]{

\begin{pspicture}(-2.5,-2)(3,2)\footnotesize \label{figure: simple form(a)}
\pspolygon[fillstyle=solid,fillcolor=lightgray,linestyle=none,linearc=0.25](1,-2)(-1,2)(3,2)(3,-2)(1,-2)
\psline[linewidth=1.5pt](0.9,-1.8)(-0.85,1.7)
\uput[dl](-0.85,1.7){$\mu(x)+\mathfrak{t}_x^\circ$}
\uput[dl](0,0){$\mu(x)$} \rput{26.568}{
\psline(0.2,0)(0.2,0.2)(0,0.2)
\psline[arrowsize=5pt]{->}(0,0)(1,0) \uput[r]{*0}(1,0){$A^* \mu_V
(V_2)$} }
\end{pspicture}

}

 \subfigure[]{

\begin{pspicture}(-2.5,-2)(3,2)\footnotesize \label{figure: simple form(b)}
\pspolygon[fillstyle=solid,fillcolor=lightgray,linestyle=none,linearc=0.25](1,-2)(-1,2)(3,2)(3,-2)(1,-2)
\psline[linewidth=1.5pt](0.9,-1.8)(-0.85,1.7)
\uput[dl](-0.85,1.7){$\mu(x)+\mathfrak{t}_x^\circ$}
\uput[dl](0,0){$\mu(x)$}

\rput{26.568}{ \psline(0.2,0)(0.2,0.2)(0,0.2)
\psline[arrowsize=5pt]{->}(0,0)(1,0)
\uput[dl]{*0}(1,0){$A^* \mu_V (V_2)$}
\psline[arrowsize=5pt]{->}(0,0)(2,0)
\uput[u]{*0}(2,0){$A^* \mu_V (V_1)$} }
\end{pspicture}

}

}

\subfigure[]{

\begin{pspicture}(-3,-2)(3,2)\footnotesize \label{figure: simple form(c)}
\pspolygon[fillstyle=solid,fillcolor=lightgray,linestyle=none,linearc=0.25](-3,-2)(3,-2)(3,2)(-3,2)
\psline[linewidth=1.5pt](1,-2)(-1,2)
\uput[dl](-1,2){$\mu(x)+\mathfrak{t}_x^\circ$}
\uput[l](0,0.1){$\mu(x)$}

\rput{26.568}{ \psline(0.2,0)(0.2,0.2)(0,0.2)
\psline[arrowsize=5pt]{->}(0,0)(1,0)
\uput[r]{*0}(1,0){$A^* \mu_V (V_2)$}

\psline(-0.2,0)(-0.2,-0.2)(0,-0.2)
\psline[arrowsize=5pt]{->}(0,0)(-1.5,0)

\uput[d]{*0}(-1.5,0){$A^* \mu_V (V_1)$} }
\end{pspicture}

}

\end{center}
\caption{\label{figure: simple form} Generalized moment maps of
Example~\ref{example: simple form of moment}}
\end{figure}

\begin{example} \label{example: simple form of moment}
We visualize an open $(n-2)$-face $\mu(M_i).$ Assume that a point
$x$ in $M_i$ has the one-dimensional stabilizer $T_x^{n-1}.$
First, let us calculate stabilizers of points in $E.$ It is easy
to show that each point $[t,\eta,0]$ in $E$ has the stabilizer
$T_x^{n-1}.$ But, to calculate stabilizers of other points we need
to know the $T_x^{n-1}$ representation $V.$ Note that the
annihilator $\mathfrak{t}_x^\circ$ is $(n-2)$-dimensional, and
hence the representation $V$ is complex two-dimensional. Since
$T_x^{n-1}$ is abelian, the representation $V$ can be written as a
sum of two complex one-dimensional subrepresentations
$$V=V_1 \oplus V_2.$$
 One of the following four cases holds. Let $H$ be the identity component of
$T_x^{n-1}.$ Note that $\mu([t,\eta,0])=\mu(x)+\eta$ and the
representation $V$ is nontrivial because the $T^{n-1}$ action is
effective.

\begin{enumerate}
    \item[i.] The subrepresentation $V_1$ is trivial.
    Then, the set of points with the stabilizer
    $T_x^{n-1}$ is the set $\{[t,\eta,v] \in E | v \in V_1 \}$
    because the $T_x^{n-1}$ representation $V_2$ is faithful.
The image of the set under $\mu$ is contained in
$\mu(x)+\mathfrak{t}_x^\circ.$ The generalized moment map image is
drawn in Figure~\ref{figure: simple form(a)}.
    \item[ii.] The subrepresentatin $V_1$ is nontrivial but a trivial $H$
    representation. Then, the set of points with the stabilizer
    $T_x^{n-1}$ is smaller than the set $\{[t,\eta,v] \in E | v \in V_1
    \}.$ The image $\mu([t,\eta,v])$ for $v \in V_1$ is also contained in $\mu(x)+\mathfrak{t}_x^\circ.$
    Also, see Figure~\ref{figure: simple form(a)}.
    \item[iii.] Both $V_1$ and $V_2$ are nontrivial $H$
    representations and $\mu_V (V) \ne \mathfrak{t}^*_x.$ Then, the set of points with the stabilizer
    $T_x^{n-1}$ is the set $\{ [t,\eta,0] \in E \}.$ The generalized moment map image is
drawn in Figure~\ref{figure: simple form(b)}.
    \item[iv.] Both $V_1$ and $V_2$ are nontrivial $H$
    representations and $\mu_V (V) = \mathfrak{t}^*_x.$ Then, the set of points with the stabilizer
    $T_x^{n-1}$ is the set $\{ [t,\eta,0] \in E \}.$ The generalized moment map image is
drawn in Figure~\ref{figure: simple form(c)}.
\end{enumerate}
Note that in all four cases the open $(n-2)$-face $\mu(M_i)$ is
locally of the form $\mu(y)+\mathfrak{t}_x^\circ$ for each $y \in
M_i.$
\end{example}

Now, we investigate local behavior of a generalized moment map
near a fixed point through describing the x-ray. To do it, we need
to calculate stabilizers near a fixed point. The $T^{n-1}$
representation $\T_x M$ for $x \in M^{T^{n-1}}$ can be expressed
as $$(\alpha_1 z_1, \cdots, \alpha_n z_n)$$ and the moment map
$\mu$ is written as $$ \mu(z_1, \cdots, z_n) = \sum \frac 1 2
|z_i|^2 \alpha_i. $$ where $\alpha_i$'s are the weights of $\T_x
M$ and $z_i \in \C.$ Then, we can easily verify that $\mu (\T_x M)
= S(\alpha_1, \cdots, \alpha_n).$ For a while, we assume that the
cone $S(\alpha_1, \cdots, \alpha_n)$ is equal to ${\mathfrak
t}^*.$ From this condition, any $(n-1)$ weights of $\alpha_1,
\cdots, \alpha_n$ are independent in ${\mathfrak t}^*.$ In other
words, any $(n-1)$-dimensional matrix subrepresentation has a
finite kernel. Also, we can regard an $(n-1)$-dimensional matrix
subrepresentation as a homomorphism from $T^{n-1}$ to
${T^{n-1}}^\prime$ where the latter group ${T^{n-1}}^\prime $ acts
diagonally on $\C^{n-1}$ by $(t_1 z_1 , \cdots , t_{n-1} z_{n-1})$
for $t_i \in S^1.$ From this, we obtain the following lemma. For
 $1 \le i < j \le n,$ put $$W_{i,j} = \{ (z_1, \cdots, z_n)\in
\C^n | z_i =z_j =0, z_k \ne 0 {\rm{~for~}} k \ne i , j \}.$$

\begin{lemma} \label{lemma: one-dimensional stabilizer}
For $n \ge 3,$ any point in $W_{i, j}$ has the same stabilizer.
The set of points in $\T_x M$ whose stabilizer is one-dimensional
is the union of $W_{i,j}$ for $1 \le i < j \le n.$
\end{lemma}

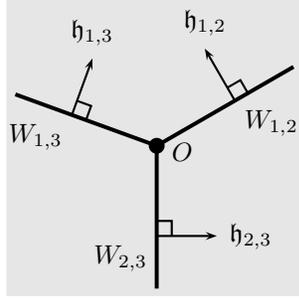
\begin{figure}[ht]
\begin{center}
\begin{pspicture}(-2,-2)(2,2)\footnotesize
\pspolygon[fillstyle=solid,fillcolor=lightgray,linestyle=none](-2,-2)(2,-2)(2,2)(-2,2)

\rput{30}{ \psline[linewidth=1.5pt](0,0)(2.1,0)
\psline{->}(1.2,0)(1.2,0.8) \psline(1.4,0)(1.4,0.2)(1.2,0.2)
\uput[dl]{*0}(1.2,0){$W_{1,2}$}
\uput[u]{*0}(1.2,0.8){$\mathfrak{h}_{1,2}$} }

\rput{160}{ \psline[linewidth=1.5pt](0,0)(2,0)
\psline{->}(1.2,0)(1.2,-0.8) \psline(1,0)(1,-0.2)(1.2,-0.2)
\uput[dr]{*0}(1.2,0){$W_{1,3}$}
\uput[u]{*0}(1.2,-0.8){$\mathfrak{h}_{1,3}$} }

\rput{270}{ \psline[linewidth=1.5pt](0,0)(1.9,0)
\psline{->}(1.2,0)(1.2,0.8) \psline(1,0)(1,0.2)(1.2,0.2)
\uput[dr]{*0}(1.2,0){$W_{2,3}$}
\uput[r]{*0}(1.2,0.8){$\mathfrak{h}_{2,3}$} }

\pscircle[fillstyle=solid,fillcolor=black](0,0){3pt}
\put(0.2,-0.2){$O$}
\end{pspicture}
\end{center}
\caption{\label{stabilizer}$W_{i,j}$ and ${\mathfrak h}_{i,j}$ of
Remark~\ref{remark: open 1-face near fixed point}}
\end{figure}

\begin{remark} \label{remark: open 1-face near fixed point}
Easily, we get $\mu(W_{i,j}) = S^\circ (\alpha_1 , \cdots ,
\hat{\alpha_i}, \cdots , \hat{\alpha_j} ,\cdots , \alpha_n )$
where the hat means a missing part. Let $H_{i, j}$ be the circle
subgroup of $T^{n-1}$ such that $W_{i, j} \subset (\T_x M)^{H_{i,
j}},$ i.e. each $H_{i,j}$ is the identity component of the
stabilizer of $W_{i,j}.$ Since any $(n-1)$ weights of $\alpha_1,
\cdots, \alpha_n$ are linearly independent, for $\{i,j\} \ne
\{k,l\}$ the images $\mu(W_{i,j})$ and $\mu(W_{k,l})$ span
different vector spaces, respectively. Thus, by
Example~\ref{example: simple form of moment} or \cite[Theorem
3.6]{GS}, we obtain $H_{i, j} \ne H_{k, l}$ for $\{i, j\} \ne \{k,
l\}.$ Also, see Figure~\ref{stabilizer}.
\end{remark}

\begin{remark} \label{remark: local picture of one-dimensional stabilizer}
Let $M_i$ be an element of $\mathfrak{M}$ with a one-dimensional
stabilizer $T_i,$ and $H$ be the identity component of $T_i.$
Assume that the closure $\overline{M_i}$ contains a fixed point
$x\in M^{T^{n-1}}$ such that $S(\alpha_{1,x}, \cdots,
\alpha_{n,x})={\mathfrak t}^*.$ By Lemma~\ref{lemma:
one-dimensional stabilizer}, $\dim M_i = 2(n-2).$ Forgetting the
whole group $T^{n-1}$ action on $M^{2n},$ we only consider the $H$
action on $M^{2n}.$ Since $H$ fixes $\overline{M_i},$ $(n-2)$
weights of $\alpha_{1,x}|_{\mathfrak h}, \cdots,
\alpha_{n,x}|_{\mathfrak h}$ are zero as elements of ${\mathfrak
h}^*.$ Since $S(\alpha_{1,x}, \cdots, \alpha_{n,x})={\mathfrak
t}^*$ implies $S(\alpha_{1,x}|_{\mathfrak h}, \cdots,
\alpha_{n,x}|_{\mathfrak h})={\mathfrak h}^*,$ the remaining two
nonzero weights have different signs in ${\mathfrak h}^*.$
Therefore, the $H$ representation on the fiber at $x$ of the
normal bundle of $\overline{M_i}$ can be expressed as $(t^d z_1,
t^{-d^\prime} z_2 )$ for $t \in H$ where $d$ and $d^\prime$ are
positive integers. Also, since $H$ fixes $\overline{M_i},$ the $H$
representation on each fiber of the normal bundle is all the same,
and we can show that it is $V$ in Example~\ref{example: simple
form of moment}. This gives a description of a generalized moment
map near $M_i$ by Example~\ref{example: simple form of moment}-iv.
Moreover, these arguments also applies to the case
$S(\alpha_{1,x}, \cdots, \alpha_{n,x}) \ne {\mathfrak t}^*.$ From
this, one may say that the cone $S(\alpha_{1,x}, \cdots,
\alpha_{n,x})$ of a fixed point $x$ in $\overline{M_i}$ determines
the generalized moment map image near $\overline{M_i}.$
\end{remark}

Next, we investigate the image(not the x-ray) of $\mu$ for the
$T^{n-1}$ representation $\T_x M$ when $x \in M^{T^{n-1}}$ and
$S(\alpha_1, \cdots, \alpha_n) \ne {\mathfrak t}^*.$ For $v \in
\mathfrak{t},$ we denote by $v^\circ$ the annihilator of $v$ in
$\mathfrak{t}^*$ with respect to the pairing between $\mathfrak t$
and ${\mathfrak t}^*$.

\begin{lemma} \label{lemma: corners in lie algebra}
Let $\alpha_1, \cdots, \alpha_n$ be the weights of a $T^{n-1}$
representation on $\C^n.$ Assume that $S(\alpha_1, \cdots,
\alpha_n) \ne {\mathfrak t}^*.$ Then, the cone $S(\alpha_1,
\cdots, \alpha_n)$ satisfies either
\begin{enumerate}
    \item[i.] it contains no nontrivial vector
    space, and there is a vector $v_0$ in $\Lambda_0$ such
    that $v_0^\circ \cap S(\alpha_1, \cdots, \alpha_n)=0.$
    \item[\textit{or} ii.] it contains a nontrivial vector
    space, and the cone is homeomorphic to $\R^l \times
(\R^1_{\ge 0})^{n-1-l}$ for some $l \in \N.$ Also, there is a
vector $v_0$ in $\Lambda_0$ such that $v_0^\circ \cap S(\alpha_1,
\cdots, \alpha_n)$ is an $l$-dimensional vector space.
\end{enumerate}

\end{lemma}

\begin{proof}
For a point $x$ having the identity as $T_x^{n-1},$ the moment map
is a submersion near $x$ by the simple form of $\mu.$ Therefore,
weights $\alpha_1, \cdots, \alpha_n$ span ${\mathfrak t}^*$ as a
vector space because the moment map image is equal to $S(\alpha_1,
\cdots, \alpha_n).$ From this, we can find $(n-1)$ independent
weights(say $\alpha_1, \cdots, \alpha_{n-1}$). Put
$\alpha_n=\sum_{i=1}^{n-1} d_i \alpha_i.$ By a suitable linear
transformation, we may assume that each $\alpha_i$ is the usual
basis element $e_i$ of $\R^{n-1}$ for $i=1, \cdots, n-1.$ In the
coordinates, $\alpha_n=(d_1, \cdots, d_{n-1}).$

Now, assume that $S(\alpha_1, \cdots, \alpha_n)$ contains no
nontrivial vector space. Then, $\alpha_n$ is the zero vector or at
least one of $d_i$'s is positive. If $\alpha_n$ is the zero
vector, then take $v_0=(-1, \cdots ,-1)$ and we are done.
Otherwise, we can find a vector $v_0=(a_1, \cdots, a_{n-1})$ with
negative integers $a_i$ such that all the pairings $\langle v_0 ,
\alpha_i \rangle$ are negative. Then, $v_0^\circ \cap S(\alpha_1,
\cdots, \alpha_n)=0.$ If $v_0$ is not contained in $\Lambda_0,$
then we only have to multiply $v_0$ by a big natural number to
obtain a wanted vector in $\Lambda_0.$

Next, assume that $S(\alpha_1, \cdots, \alpha_n)$ contains a
nontrivial vector space. This implies that each $d_i$ is
nonpositive. Let $l$ be the number of negative $d_i$'s(say $d_1,
\cdots, d_l$). Then, $S(\alpha_1, \cdots, \alpha_n)$ is equal to
$\R^l \times (\R^1_{\ge 0})^{n-1-l}.$ Put $v_0=(a_1, \cdots,
a_{n-1})$ with $a_1, \cdots, a_l=0$ and $a_{l+1}, \cdots,
a_{n-1}=-1.$ Then, $v_0^\circ \cap S(\alpha_1, \cdots,
\alpha_n)=\R^l \times 0.$

\end{proof}

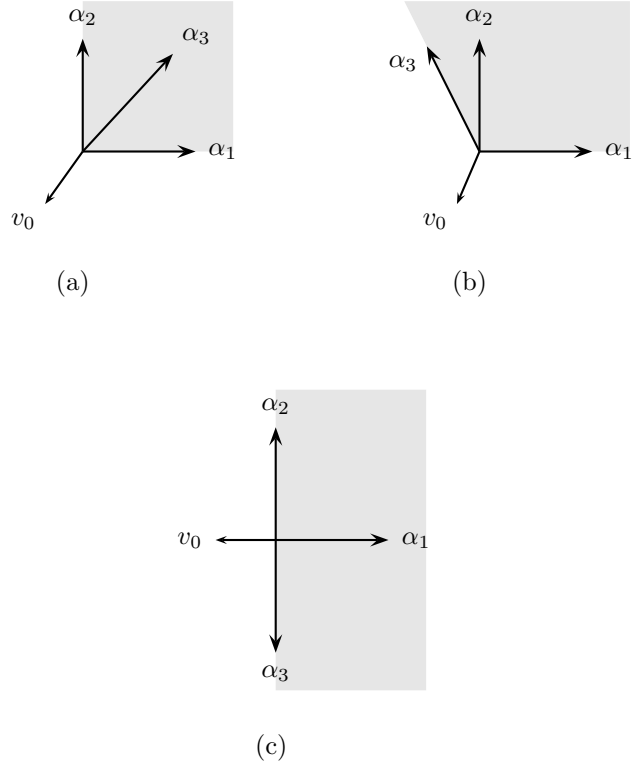
\begin{figure}[ht]
\begin{center}

\mbox{
\subfigure[]{
\begin{pspicture}(-2.5,-1)(2.5,2)\footnotesize \label{figure: three cases of the cone(a)}
\pspolygon[fillstyle=solid,fillcolor=lightgray,linestyle=none](0,0)(2,0)(2,2)(0,2)
\psline[arrowsize=5pt]{->}(0,0)(1.5,0) \uput[r](1.5,0){$\alpha_1$}
\psline[arrowsize=5pt]{->}(0,0)(0,1.5) \uput[u](0,1.5){$\alpha_2$}
\psline[arrowsize=5pt]{->}(0,0)(1.2,1.3)
\uput[ur](1.2,1.3){$\alpha_3$} \psline{->}(0,0)(-0.5,-0.7)
\uput[dl](-0.5,-0.7){$v_0$}

\end{pspicture}}

\subfigure[]{
\begin{pspicture}(-2.5,-1)(2.5,2)\footnotesize \label{figure: three cases of the cone(b)}
\pspolygon[fillstyle=solid,fillcolor=lightgray,linestyle=none](0,0)(2,0)(2,2)(-1,2)
\psline[arrowsize=5pt]{->}(0,0)(1.5,0) \uput[r](1.5,0){$\alpha_1$}
\psline[arrowsize=5pt]{->}(0,0)(0,1.5) \uput[u](0,1.5){$\alpha_2$}
\psline[arrowsize=5pt]{->}(0,0)(-0.7,1.4)
\uput[dl](-0.7,1.4){$\alpha_3$} \psline{->}(0,0)(-0.3,-0.7)
\uput[dl](-0.3,-0.7){$v_0$}

\end{pspicture}}

}

\subfigure[]{\begin{pspicture}(-2.5,-2)(2.5,2.5)\footnotesize
\label{figure: three cases of the cone(c)}
\pspolygon[fillstyle=solid,fillcolor=lightgray,linestyle=none](0,-2)(2,-2)(2,2)(0,2)
\psline[arrowsize=5pt]{->}(0,0)(1.5,0) \uput[r](1.5,0){$\alpha_1$}
\psline[arrowsize=5pt]{->}(0,0)(0,1.5) \uput[u](0,1.5){$\alpha_2$}
\psline[arrowsize=5pt]{->}(0,0)(0,-1.5)
\uput[d](0,-1.5){$\alpha_3$} \psline{->}(0,0)(-0.8,0)
\uput[l](-0.8,0){$v_0$}

\end{pspicture}}

\end{center}
\caption{\label{figure: three cases of the cone} Three cases of
$S(\alpha_1, \alpha_2, \alpha_3)$ when $\alpha_3 \ne 0$}
\end{figure}

\begin{example} 
Figure~\ref{figure: three cases of the cone} illustrates
Lemma\ref{lemma: corners in lie algebra} when $n=3$ and $\alpha_3
\ne 0.$ The case of $\alpha_3=0$ is the same with
Figure~\ref{figure: three cases of the cone(a)} except
$\alpha_3=0.$
\end{example}

\section{Proof} \label{section: 3}

First, we prove the following lemma about an equivalent condition
to guarantee that a symplectic action must be Hamiltonian.

\begin{lemma} \label{lemma: equivalent condition}
Let an $m$-dimensional torus $T^m$ act on $(M^{2n}, \omega).$ Then
the action is Hamiltonian if and only if there exists a component
$F$ of $M^{T^m}$ such that the map $\pi_1 (F) \rightarrow \pi_1
(M^{2n})$ is surjective.
\end{lemma}

\begin{proof}
Assume that the action is Hamiltonian. Let $F$ be $\mu^{-1}(p)$
where $\mu$ is a moment map for the action and the point $p$ is an
extreme point(or a vertex) of the image of $\mu.$ For an extreme
point, see \cite[Chapter 3]{R}. Then, $F$ is a component of
$M^{T^m}$ and connected by Convexity Theorem \cite{A}, \cite{GS}.
By \cite[p. 14 (13)]{Fu}, there exists a vector $v_0 \in \Lambda_0
\subset \mathfrak t$ such that $(p + v_0^\circ ) \cap
\mu(M^{2n})=p$ since $p$ is an extreme point of $\mu(M^{2n}).$ Let
$H$ be the circle subgroup of $T^m$ corresponding to $v_0.$ Then,
a moment map $\mu_H$ for the $H$ action attains its extremum at
$F,$ and hence $F$ is also a component of $M^H.$ By \cite{L}, the
fundamental groups $\pi_1 (M^{2n}), \pi_1 (M_{min}), \pi_1
(M_{max})$ are all the same where $M_{min}$ and $M_{max}$ mean the
preimages of the minimum and maximum under $\mu_H,$ respectively.
Thus, $\pi_1 (F) \rightarrow \pi_1 (M^{2n})$ is surjective.

We prove the converse. Assume that for a circle subgroup $H$ of
$T^m,$ the $H$ action is non-Hamiltonian. Then, there exists a
loop in $M^{2n}$ such that its image under a generalized moment
map $\mu_H$ for the $H$ action is nontrivial in $\pi_1(\R/\Z).$
But, the loop can not be deformed into a component $F$ of
$M^{T^m}$ because $\mu_H (F)$ is a point in $\R/\Z.$
\end{proof}

We prove the theorem by induction on $n.$ The case of $n=2$ is
McDuff's theorem \cite{M}. Suppose that the theorem is true when
the dimension of the manifold $M$ is less than $2n.$ The situation
consists of the following two cases.


\begin{figure}[ht]
\begin{center}
\begin{pspicture}(-2,-2)(2,1)\footnotesize
\pspolygon[fillstyle=solid,fillcolor=lightgray,linestyle=none](-2,-2)(2,-2)(2,0)(-2,0)
\uput[ur](-2,0){$N$}

\psline[linewidth=1.5pt](-2,0)(2,0)

\psline[arrowsize=6pt]{->}(0,0)(1.2,0) \uput[d](1.2,0){$\alpha_3$}

\psline[arrowsize=6pt]{->}(0,0)(-1.2,0)
\uput[d](-1.2,0){$\alpha_1$}

\psline[arrowsize=6pt]{->}(0,0)(-0.8,-1.2)
\uput[d](-0.8,-1.2){$\alpha_2$}

\psline{->}(0,0)(0,0.8) \uput[u](0,0.8){$v_0$}

\psline(0.2,0)(0.2,0.2)(0,0.2)

\end{pspicture}
\end{center}
\caption{\label{case1} Proof of Case 1}
\end{figure}
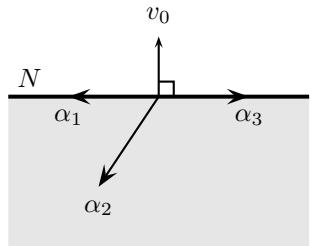

\emph{Case 1:} There exists a fixed point $x$ such that
$S(\alpha_{1,x}, \cdots, \alpha_{n,x}) \ne {\mathfrak t}^*.$

By Lemma~\ref{lemma: corners in lie algebra}, there exists a
vector $v_0 \in \Lambda_0 \subset \mathfrak t$ such that the
pairing $\langle v_0, S(\alpha_{1,x}, \cdots, \alpha_{n,x})
\rangle$ is nonpositive and $v_0^\circ \cap S(\alpha_{1,x},
\cdots, \alpha_{n,x}),$ say $Z,$ is an $l$-dimensional vector
space for some $0 \le l \le n-2.$ We denote by $H$ the circle
subgroup of $T^{n-1}$ corresponding to $v_0.$ Then, a generalized
moment map $\mu_H$ of the $H$ action attains its local maximum,
and hence the $H$ action is Hamiltonian as noted in the proof of
\cite[Lemma 2]{M}. So, we can assume that $\mu_H$ is an $\R$
valued moment map with the maximum $0.$ If we put
$N=\mu_H^{-1}(0),$ then $\pi_1 (N) \rightarrow \pi_1 (M^{2n})$ is
surjective by \cite{L}. If $l=0$ and all $\alpha_{i,x}$ are
nonzero, then $N=\{x\}$ and we obtain the proof by
Lemma~\ref{lemma: equivalent condition}. Similarly, if $l=0$ and
one of $\alpha_{i,x}$'s is zero, then $N$ is a two-surface fixed
by $T^{n-1}$ and the proof is obtained in the same way.

Now, assume that $l \ge 1$. Then, $\dim N=2(l+1)$ by the proof of
Lemma~\ref{lemma: corners in lie algebra}. Also, $N$ is fixed by
the $(n-1-l)$-dimensional subtorus $exp(Z^\circ)$ generated by
$Z^{\circ}.$ There is an $l$-dimensional subtorus of $T^{n-1}$
denoted by $T^l$ whose Lie algebra and $Z^{\circ}$ span $\mathfrak
t,$ i.e. $T^{n-1}=T^l \cdot exp(Z^\circ).$ The subtorus $T^l$ acts
symplectically on $N$ with a nonempty fixed point set containing
$x$ because $x \in M^{T^{n-1}}.$ Thus, the $T^l$ action on $N$ is
Hamiltonian by the induction hypothesis. Therefore, any loop in
$N$ can be deformed into a component $F$ of $N^{T^l}=N^{T^{n-1}},$
i.e. $\pi_1 (F) \rightarrow \pi_1 (N)$ is surjective. Since we
already have the surjection $\pi_1 (N) \rightarrow \pi_1
(M^{2n}),$ the proof is obtained again by Lemma~\ref{lemma:
equivalent condition}.


\emph{Case 2:} For each fixed point $x,$ we have $S(\alpha_{1,x},
\cdots, \alpha_{n,x})=\mathfrak t^*.$

To obtain the proof, we need to prove the following lemma. Let
$M_i$ be an element of $\mathfrak M$ with a one-dimensional
stabilizer $T_i.$ We prove that if the generalized moment map
$\mu$ behaves like one of Figure~\ref{figure: simple form(a)},
\ref{figure: simple form(b)} near a point in $M_i,$ then the
closure $\overline{M_i}$ contains a point in $M^{T^{n-1}}.$ But,
this can not happen under the assumption of Case 2 by
Remark~\ref{remark: local picture of one-dimensional stabilizer}.

\begin{figure}[ht]
\begin{center}
\begin{pspicture}(-2,-2)(2,3)\footnotesize
\pspolygon[fillstyle=solid,fillcolor=lightgray,linestyle=none](-2,-2)(2,-2)(2,2)(-2,2)
\psline[linewidth=1.5pt](-2,2)(2,2)
\uput[ur](-2,2){$\overline{M_i}$}

\psline{->}(-2.5,-2)(-2.5,2.5) \uput[l](-2.5,2){$\mu_H(x)$}
\psdots(-2.5,2)(-2.5,-1)(-2.5,-0.2)(-2.5,0.6)
\uput[l](-2.5,2){$\mu_H(x)$}

\uput[l](-2.5,-1){$\mu_H(x_1)$}

\uput[l](-2.5,-0.2){$\mu_H(x_2)$}

\uput[l](-2.5,0.6){$\mu_H(x_3)$}

\uput[u](-2.5,2.5){$\mu_H$}

\psline[linewidth=1.5pt](-0.8,-2)(-0.5,-1)(-1.2,-0.5)

\psline[linewidth=1.5pt](-0.5,-1)(0.1,-0.2)(1,-1.2)

\psline[linewidth=1.5pt](0.1,-0.2)(-0.2,0.6)(-1,1)

\psline[linewidth=1.5pt](-0.2,0.6)(0,1)

\psline[linewidth=1.5pt,linestyle=dotted,dotsep=1pt](-0.2,1.2)(-0.2,1.8)

\psdots[dotsize=5pt](-0.5,-1)(0.1,-0.2)(-0.2,0.6)

\uput[dr](-0.5,-1){$x_1$}

\put(-0.7,-0.5){$W_1$}

\uput[ur](0.1,-0.2){$x_2$}

\uput[r](-0.2,0.6){$x_3$}

\end{pspicture}
\end{center}
\caption{\label{proof}Proof of Lemma~\ref{lemma: impossible
cases}}
\end{figure}
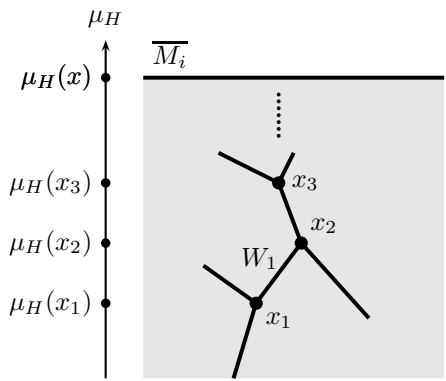

\begin{lemma} \label{lemma: impossible cases}
Let $x$ be a point in $M_i$ with a one-dimensional stabilizer
$T_i.$ If $\mu_V(V) \ne \mathfrak{t}^*_x,$ then $\overline{M_i}$
contains a point in $M^{T^{n-1}}$ where $V$ is the symplectic
slice representation at $x$ and $\mu_V$ is its moment map.
\end{lemma}

\begin{proof}
Assume that $\overline{M_i}$ does not contain any point in
$M^{T^{n-1}}.$ Let $H$ be the identity component of $T_i =
T^{n-1}_x.$ Let $\mu_H$ be the generalized moment map for the $H$
action given by
$$\mu_H ([t,\eta,v])=\mu(x)+\mu_V (v)$$ near the orbit $T^{n-1}
\cdot x.$ Since $H$ fixes $\overline{M_i}$ and $\mu_V (V) \ne
{\mathfrak t}^*_x,$ the map $\mu_H$ has a local extremum at
$\overline{M_i},$ say maximum $\mu_H (x),$ and hence the $H$
action is Hamiltonian. Therefore, we can assume that the function
$\mu_H$ is $\R$ valued.

Let $x_1 \in M^{T^{n-1}}.$ Then, by the assumption  the weights of
$\T_{x_1} M$ satisfy $S(\alpha_{1,x_1}, \cdots, \alpha_{n,
x_1})={\mathfrak t}^*,$ and this implies $\mu_H (x_1) < \mu_H (x)$
because $x_1$ is not a local extremum for $\mu_H.$ Since
$S(\alpha_{1,x_1}, \cdots, \alpha_{n, x_1})={\mathfrak t}^*,$
there is a $W_{k,l}$ of Lemma~\ref{lemma: one-dimensional
stabilizer} in $\T_{x_1} M$ such that $x_1 \in \overline{W}_{k,l}$
and $\mu_H(x_1) < \mu_H(y)$ for some $y \in \overline{W}_{k,l}.$
Let $T_{i_1}$ be the stabilizer of $W_{k,l},$ and $W_1$ be the
component of $M_{T_{i_1}}$ such that $\overline{W}_1$ contains
$x_1,$ i.e. $\overline{W}_{k,l}$ is a linear approximation of
$\overline{W}_1$ near $x_1.$ Then, the induced $T^{n-1}/T_{i_1}$
action on $\overline{W}_1$ is effective and Hamiltonian by the
induction hypothesis. Also, $\overline{W}_1$ has another fixed
point $x_2$ for the $T^{n-1}$ action such that $\mu_H(x_1) <
\mu_H(x_2)$ because $\mu_H(x_1)$ is not the maximum of
$\mu_H(\overline{W}_1).$ Since $\mu_H (x_2) < \mu_H (x),$ we can
repeat this process infinite times to obtain a sequence
$\{x_i\}_{i \in \N} \subset M^{T^{n-1}}$ such that
$$ \mu_H(x_1) < \mu_H(x_2) < \mu_H(x_3) < \cdots.$$
This is a contradiction.
\end{proof}

\begin{figure}[ht]
\begin{center}
\begin{pspicture}(-3,-2.2)(3,2)\footnotesize
\pspolygon[fillstyle=solid,fillcolor=lightgray,linestyle=dotted,dotsep=2pt,
linewidth=0.5pt](-2,-2)(2,-2)(2,2)(-2,2)

\psline{->}(-2,-2.5)(2,-2.5)

\uput[r](2,-2.5){$\mu_1$}

\psline{->}(-2.5,-2)(-2.5,2)

\uput[u](-2.5,2){$\mu_2$}

\pspolygon[linewidth=1.5pt](1.2,.8)(.8,1.2)(-.8,1.2)(-1.2,.8)(-1.2,-.8)(-.8,-1.2)(.8,-1.2)(1.2,-.8)
\psline[linewidth=1.5pt](1.2,0.8)(2.0,0.8)

\psline[linewidth=1.5pt](0.8,1.2)(0.8,2)

\psline[linewidth=1.5pt](-0.8,1.2)(-0.8,2)

\psline[linewidth=1.5pt](-1.2,0.8)(-2,0.8)

\psline[linewidth=1.5pt](-1.2,-0.8)(-2,-0.8)

\psline[linewidth=1.5pt](-0.8,-1.2)(-0.8,-2)

\psline[linewidth=1.5pt](0.8,-1.2)(0.8,-2)

\psline[linewidth=1.5pt](1.2,-0.8)(2,-0.8)

\psdots[dotsize=5pt](1.2,.8)(.8,1.2)(-.8,1.2)(-1.2,.8)(-1.2,-.8)(-.8,-1.2)(.8,-1.2)(1.2,-.8)

\psline[linewidth=1.5pt, linestyle=dashed](-2,0.25)(2,0.25)

\psdot(-2.5,0.25) \uput[l](-2.5,0.25){$\xi_2$}

\pscircle(-1.2,0.25){4pt}

\uput[dr](-1.2,0.25){$\mu(x_1)$}

\pscircle(1.2,0.25){4pt}

\uput[dl](1.2,0.25){$\mu(x_2)$}

\end{pspicture}
\end{center}
\caption{\label{the idea}Illustration of the idea of proof}
\end{figure}
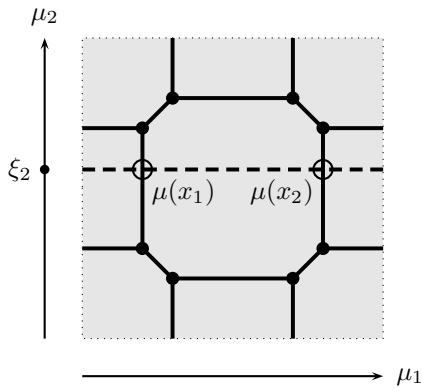

Before we start the proof of Case 2, we will illustrate the idea
of proof with the following example of the case $n=3.$ Assume that
there exists a six-dimensional two-torus action whose x-ray is
Figure~\ref{the idea}. Note that $S(\alpha_{1,x}, \alpha_{2,x},
\alpha_{3,x})={\mathfrak t}^*$ for each $x \in M^{T^2}.$ Let us
fix a decomposition of the two-torus $T^2=S^1_1 \times S^1_2,$ and
let $\R/\Z$ valued functions $\mu_i$ be generalized moment maps
for the $S^1_i$ actions. Put $\mu=(\mu_1, \mu_2).$ The gray square
is the range $\R/\Z \times \R/\Z$ of $\mu.$ Each black dot is the
image of a fixed point, and each thick line is the image of an
open one-face. Let $\xi_2 \in \R/\Z$ be a regular value for
$\mu_2.$ Then, the thick dashed line becomes the x-ray of the
induced symplectic $T^2/S^1_2$ action on the four-dimensional
orbifold $\mu^{-1}_2 (\xi_2)/S^1_2.$ The induced action has two
isolated fixed points which are orbits of $x_1$ and $x_2.$ Also,
the induced action has the generalized moment map induced from
$\mu_1,$ and the induced generalized moment map can not be an $\R$
valued function because two fixed points are not locally extremal.
But, it is impossible by the same proof of \cite[Lemma 3]{M} in
the orbifold setting. Thus, there is no action with
Figure~\ref{the idea} as an x-ray.

Now, we start the proof. We can choose a generic regular value
$\xi_i \in \R/\Z$ for each $\mu_i,$ $i=2, \cdots, n-1,$ such that
in the range of $\mu,$ the subset $(\R/\Z) \times \xi_2 \times
\cdots \times \xi_{n-1}$ meets each open $m$-face transversely for
$m=1, \cdots, n-2$ and there is at least one open $(n-2)$-face
intersecting the subset. Transversality guarantees that the subset
$(\R/\Z) \times \xi_2 \times \cdots \times \xi_{n-1}$ does not
meet any open $m$-face for $m \le n-3,$ and this is why our
interest is mainly in (open) $(n-2)$-faces.

Put $\xi=(\xi_2, \cdots, \xi_{n-1}).$ Also, we put $V_\xi=( \mu_2,
\cdots, \mu_{n-1} )^{-1}(\xi).$ Then, the orbit space $B_\xi=V_\xi
/ (S_2^1 \times \cdots \times S_{n-1}^1)$ is a four-dimensional
orbifold naturally endowed with a symplectic form. For orbifolds
see \cite{Sa}, and for symplectic actions on orbifolds see
\cite{LT}. Denote $T^{n-1}/(S_2^1 \times \cdots \times S_{n-1}^1)$
by $H.$ If we restrict the generalized moment map $\mu_1$ to
$V_\xi,$ then we obtain the induced generalized moment map $\mu_H
: B_\xi \rightarrow \R/\Z$ for the induced symplectic $H$ action
on $B_\xi.$ Also, the $H$ action on $B_\xi$ is effective for a
generic $\xi.$ However, we do not know whether $\mu_H$ can be
regarded as an $\R$ valued function.

Let $\pi: V_\xi \rightarrow B_\xi$ be the orbit map. For $x \in
V_\xi,$ the point $\pi(x)$ is an $H$ fixed point if and only if
the orbit $T^{n-1} \cdot x$ is $(n-2)$-dimensional because $S_2^1
\times \cdots \times S_{n-1}^1$ acts almost freely on $V_\xi,$
i.e. stabilizers are all finite groups, and hence $H \cdot x
\subset (S_2^1 \times \cdots \times S_n^1) \cdot x.$ This is also
equivalent to say that the stabilizer of $x$ is one dimensional.
Since $(\R/\Z) \times \xi$ meets at least one open $(n-2)$-face
and an open $(n-2)$-face is an image of some $M_i \in \mathfrak M$
with a one-dimensional $T_i,$ $B_\xi$ has at least one fixed
point.

\begin{figure}[ht]
\begin{center}

\begin{pspicture}(-3,-2)(3,2)\footnotesize
\pspolygon[fillstyle=solid,fillcolor=lightgray,linestyle=none,linearc=0.25](-2,-2)(2,-2)(2,2)(-2,2)
\psline[linewidth=1.5pt](1,-2)(-1,2)
\uput[r](-0.8,1.8){$\mu(x)+\mathfrak{t}_x^\circ$}

\psdot(0,0) \uput[ur](0,0){$\mu(x)$}

\psline[linewidth=1.5pt, linestyle=dashed](-2,0)(2,0)
\uput[r](2,0){$(\R/\Z) \times \xi$}

\end{pspicture}

\end{center}
\caption{\label{figure: proof of case 2} Proof of Case 2}
\end{figure}
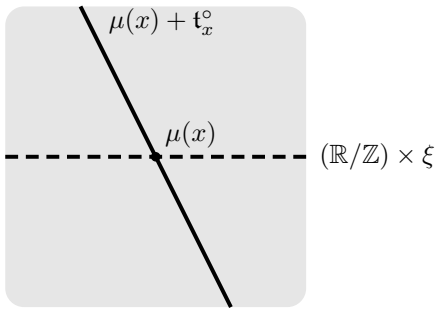

Let $\pi(x) \in B_\xi$ be a fixed point of the $H$ action. The
orbit space $B_\xi$ need not be connected because $V_\xi$ need not
be connected. So, we denote the component of $B_\xi$ containing
$\pi(x)$ by the same notation $B_\xi.$ From Remark~\ref{remark:
local picture of one-dimensional stabilizer} and Lemma~\ref{lemma:
impossible cases}, the x-ray of the $T^{n-1}$ action on $M^{2n}$
is locally of the form of Figure~\ref{figure: simple form(c)} near
$\mu(x).$ Also, see Figure~\ref{figure: proof of case 2}.
Therefore, we can conclude that each fixed point in $B_\xi$ is an
isolated critical point because it is not extremal for $\mu_H.$


In summary, we obtain a symplectic circle action on a
four-dimen\-sional symplectic orbifold with nonempty isolated
fixed points which are nonextremal for a generalized moment map.
McDuff shows that such an action does not exist if the orbifold is
a smooth manifold \cite[Lemma 3]{M}. And, it is easy to see that
her proof can be transferred almost literally to our case, and we
obtain a contradiction. Therefore, Case 2 can not happen.


As a corollary of the above proof, we obtain the following. We
define the dimension of the fixed point set of an action as the
maximum of dimensions of components of the fixed point set of the
action.

\begin{corollary}
Let $T^m$ act symplectically on $(M^{2n}, \omega)$ with a nonempty
$M^{T^m}$ in an effective way. If $\dim M^{T^m} \ge 2(n-m),$ then
the action is Hamiltonian. In particular, if $m=n,$ then the
action has a fixed point if and only if it is Hamiltonian.
\end{corollary}

\begin{proof}
Let $F$ be a component of the fixed point set such that $\dim
F=\dim M^{T^m}$. By the assumption, $\dim F \ge (2n-2m).$ Then at
any point $x \in F,$ there are $n-\dim F/2$ nonzero weights of
$\T_x M.$ Since $m \ge n-\dim F/2$ and the weights span
${\mathfrak t}^*$ as a vector space, the number $n-\dim F/2$ must
be equal to $m.$ The nonzero weights are linearly independent in
$\mathfrak t^*,$ and hence $S(\alpha_{1,x}, \cdots, \alpha_{n,x})$
is strictly convex, i.e. it contains no nontrivial vector space.
By the same arguments of the proof of the case 1, we obtain the
desired proof.
\end{proof}

\end{document}